\documentclass[11pt,final]{amsart}
\usepackage{amsopn}
\usepackage{amssymb, amscd}
\usepackage{showkeys}
\usepackage{multirow}
\usepackage{graphicx, graphics, epsfig}
\usepackage{faktor}
\usepackage{enumerate}
\usepackage{tikz}
\usepackage{pgfplots}
\usepackage{hyperref}
\usepackage{color}
\usepackage{amsmath}
\usepackage{cleveref}

\usepackage{hyperref} %VER COMO SE usA-

\newcommand{\nc}{\newcommand}

\nc{\fg}{\mathfrak{f} } \nc{\vg}{\mathfrak{v} } \nc{\wg}{\mathfrak{w} }
\nc{\zg}{\mathfrak{z} } \nc{\ngo}{\mathfrak{n} } \nc{\kg}{\mathfrak{k} }
\nc{\mg}{\mathfrak{m} } \nc{\bg}{\mathfrak{b} } \nc{\ggo}{\mathfrak{g} } \nc{\eg}{\mathfrak{e} }
\nc{\ggob}{\overline{\mathfrak{g}} } \nc{\sog}{\mathfrak{so} }
\nc{\sug}{\mathfrak{su} } \nc{\spg}{\mathfrak{sp} } \nc{\slg}{\mathfrak{sl} }
\nc{\glg}{\mathfrak{gl} } \nc{\cg}{\mathfrak{c} } \nc{\rg}{\mathfrak{r} }
\nc{\hg}{\mathfrak{h} } \nc{\tg}{\mathfrak{t} } \nc{\ug}{\mathfrak{u} }
\nc{\dg}{\mathfrak{d} } \nc{\ag}{\mathfrak{a} } \nc{\pg}{\mathfrak{p} }
\nc{\sg}{\mathfrak{s} } \nc{\affg}{\mathfrak{aff} } \nc{\qg}{\mathfrak{q} } \nc{\lgo}{\mathfrak{l} }

\nc{\pca}{\mathcal{P}} \nc{\nca}{\mathcal{N}} \nc{\lca}{\mathcal{L}}
\nc{\oca}{\mathcal{O}} \nc{\mca}{\mathcal{M}} \nc{\tca}{\mathcal{T}}
\nc{\aca}{\mathcal{A}} \nc{\cca}{\mathcal{C}} \nc{\gca}{\mathcal{G}}
\nc{\sca}{\mathcal{S}} \nc{\hca}{\mathcal{H}} \nc{\bca}{\mathcal{B}}
\nc{\dca}{\mathcal{D}} \nc{\eca}{\mathcal{E}} \nc{\wca}{\mathcal{W}}

\nc{\vp}{\varphi} \nc{\ddt}{\tfrac{d}{dt}} \nc{\dsdt}{\tfrac{d^2}{dt^2}} \nc{\dds}{\frac{d}{ds}}
\nc{\dpar}{\frac{\partial}{\partial t}} \nc{\im}{\mathrm{i}}
\nc{\SO}{\mathrm{SO}} \nc{\Spe}{\mathrm{Sp}} \nc{\Sl}{\mathrm{SL}}
\nc{\SU}{\mathrm{SU}} \nc{\Or}{\mathrm{O}} \nc{\U}{\mathrm{U}} \nc{\Gl}{\mathrm{GL}}
\nc{\Se}{\mathrm{S}} \nc{\Cl}{\mathrm{Cl}} \nc{\Spin}{\mathrm{Spin}}
\nc{\Pin}{\mathrm{Pin}} \nc{\G}{\mathrm{GL}_n(\RR)} \nc{\g}{\mathfrak{gl}_n(\RR)}
\nc{\RR}{\mathbb{R}} \nc{\HH}{\mathbb{H}} \nc{\CC}{\mathbb{C}} \nc{\ZZ}{\mathbb{Z}}
\nc{\FF}{\mathbb{F}} \nc{\NN}{\mathbb{N}} \nc{\QQ}{\mathbb{Q}} \nc{\PP}{\mathbb{P}} \nc{\OO}{\mathbb{O}}

\nc{\vs}{\vspace{.2cm}} \nc{\vsp}{\vspace{1cm}} \nc{\ip}{\langle\cdot,\cdot\rangle}
\nc{\ipp}{(\cdot,\cdot)} \nc{\la}{\langle} \nc{\ra}{\rangle} \nc{\unm}{\tfrac{1}{2}}
\nc{\unc}{\tfrac{1}{4}} \nc{\und}{\frac{1}{16}} \nc{\no}{\vs\noindent}
\nc{\lam}{\Lambda^2(\RR^n)^*\otimes\RR^n} \nc{\tangz}{{\rm T}^{\rm Zar}}
\nc{\nor}{{\sf n}}  \nc{\mum}{/\!\!/} \nc{\kir}{/\!\!/\!\!/}
\nc{\Ri}{\tfrac{4\Ric_{\mu}}{||\mu||^2}} \nc{\ds}{\displaystyle}
\nc{\lb}{[\cdot,\cdot]} \nc{\isn}{\tfrac{1}{||v||^2}}
\nc{\gkp}{(\ggo=\kg\oplus\pg,\ip)} \nc{\ukh}{(\ug=\kg\oplus\hg,\ip)}
\nc{\tgkp}{(\tilde{\ggo}=\kg\oplus\pg,\ip)}
\nc{\wt}{\widetilde}
\nc{\iop}{\mathtt{i}} \nc{\jop}{\mathtt{j}}
\nc{\Hk}{H_{\kil}} \nc{\gk}{g_{\kil}}

\nc{\Hess}{\operatorname{Hess}} \nc{\ad}{\operatorname{ad}}
\nc{\Ad}{\operatorname{Ad}} \nc{\rank}{\operatorname{rk}}
\nc{\Irr}{\operatorname{Irr}} \nc{\End}{\operatorname{End}}
\nc{\Aut}{\operatorname{Aut}}\nc{\Inn}{\operatorname{Inn}}
\nc{\Der}{\operatorname{Der}} \nc{\Ker}{\operatorname{Ker}}
\nc{\Iso}{\operatorname{Iso}} \nc{\Diff}{\operatorname{Diff}}
\nc{\Lie}{\operatorname{L}} \nc{\tr}{\operatorname{tr}} \nc{\dif}{\operatorname{d}}
\nc{\sen}{\operatorname{sen}} \nc{\modu}{\operatorname{mod}}
\nc{\CRic}{\operatorname{PP}} \nc{\Cric}{\operatorname{P}} \nc{\Ricci}{\operatorname{Ric}}
\nc{\sym}{\operatorname{sym}} \nc{\herm}{\operatorname{herm}} \nc{\symac}{\operatorname{sym^{ac}}}
\nc{\symc}{\operatorname{sym^{c}}} \nc{\scalar}{\operatorname{Sc}}
\nc{\grad}{\operatorname{grad}} \nc{\ricci}{\operatorname{Rc}} \nc{\kil}{\operatorname{B}} \nc{\cas}{\operatorname{C}} \nc{\lic}{\operatorname{L}}
\nc{\Nor}{\operatorname{Norm}}  \nc{\ricc}{\operatorname{Rc^{c}}}
\nc{\Ricc}{\operatorname{Ric^{c}}} \nc{\ricac}{\operatorname{Rc^{ac}}}
\nc{\Ricac}{\operatorname{Ric^{ac}}} \nc{\Riem}{\operatorname{Rm}} \nc{\Sec}{\operatorname{Sec}}
\nc{\riccig}{\operatorname{ric^{\gamma}}} \nc{\mm}{\operatorname{m}} \nc{\Mm}{\operatorname{M}}
\nc{\Le}{\operatorname{L}} \nc{\tang}{\operatorname{T}}
\nc{\level}{\operatorname{level}} \nc{\rad}{\operatorname{r}}
\nc{\abel}{\operatorname{ab}} \nc{\CH}{\operatorname{CH}} \nc{\Cone}{{\mathcal C}} \nc{\CCone}{\operatorname{CC}} \nc{\CP}{{\mathcal P}}
\nc{\mcc}{\operatorname{mcc}} \nc{\Adj}{\operatorname{Adj}}
\nc{\Order}{\operatorname{O}}  \nc{\inj}{\operatorname{inj}} \nc{\proy}{\operatorname{pr}}
\nc{\vol}{\operatorname{vol}} \nc{\Diag}{\operatorname{Dg}} \nc{\Diagg}{\operatorname{Diag}}
\nc{\Spec}{\operatorname{Spec}} \nc{\Ima}{\operatorname{Im}} \nc{\Rea}{\operatorname{Re}}
\nc{\spann}{\operatorname{span}} \nc{\Aff}{\operatorname{Aff}} \nc{\E}{\operatorname{E}} \nc{\id}{\operatorname{id}} \nc{\dete}{\operatorname{det}} \nc{\Crit}{\operatorname{Crit}} \nc{\val}{\operatorname{val}}

\theoremstyle{plain}
\newtheorem{theorem}{Theorem}[section]
\newtheorem{proposition}[theorem]{Proposition}

\theoremstyle{definition}
\newtheorem{definition}[theorem]{Definition}

\theoremstyle{remark}
\newtheorem{remark}[theorem]{Remark}

\title{Einstein metrics on the full flag F(n).}

\author{MIKHAIL R. GUZMAN} 

\address{FaMAF, Universidad Nacional de C\'ordoba and CIEM, CONICET (Argentina)}
\email{mikhail.rios.guzman@unc.edu.ar} 

\thanks{This research was partially supported by grants from Univ. Nac. de C\'ordoba and Foncyt (Argentina)} 

\date{\today}
\pgfplotsset{compat=1.18} 
\begin{document}

%\tableofcontents

\begin{abstract}
Let $ M = G/K $ be a full flag manifold. In this work, we investigate the $ G$-stability of Einstein metrics on $M$ and analyze their stability types, including coindices, for several cases. We specifically focus on $F(n) = \mathrm{SU}(n)/T$, emphasizing $n = 5$, where we identify four new Einstein metrics in addition to known ones. Stability data, including coindex and Hessian spectrum, confirms that these metrics on $F(5)$ are pairwise non-homothetic, providing new insights into the finiteness conjecture.
\end{abstract}

\maketitle

\section{Introduction}\label{intro}

In the study of homogeneous Einstein metrics, a main open problem is about finiteness up to homothecy.  More precisely, given a compact homogeneous space $M=G/K$, the question is whether the moduli space 
$$
\eca_1^G/\Aut(G/K)
$$ 
is always finite, where $\eca_1^G$ is the set of all $G$-invariant unit volume Einstein metrics on $M$ and $\Aut(G/K)\subset\Diff(M)$ is the Lie group of automorphisms of $G$ leaving $K$ invariant.  This has been conjectured to hold in the case when the isotropy representation is multiplicity-free by B\"ohm, Wang and Ziller in \cite{BWZSoloTeoremaECompacto} (note $\eca_1^G$ should itself be finite in this case since $\Aut(G/K)$ is finite), but it is also open in the general case.  

From a variational point of view, it is well known that $g\in\mca_1^G$ is Einstein if and only if $g$ is a critical point of the scalar curvature functional 
$$
\scalar:\mca_1^G\longrightarrow \RR, 
$$ 
where $\mca^G_1$ denotes the manifold of all unit volume $G$-invariant metrics on $M$.  The tangent space is given by  
$$
T_g\mca_1^G = T_g\Aut(G/K)\cdot g \oplus \tca\tca_g^G, 
$$
where $\tca\tca_g^G:=(\Ker\delta_g\cap\Ker\tr_g)^G$ is the space of all TT-tensors which are $G$-invariant (see \cite{jLauret}). 

In this light, if an Einstein metric $g$ is {\it non-degenerate} as a critical point (i.e., $\scalar''_g|_{\tca\tca_g^G}$ is non-degenerate), then $g$ is an isolated point in the moduli space $\eca_1^G/\Aut(G/K)$, which is compact by  \cite[Theorem 1.6]{BWZSoloTeoremaECompacto}.  This provides a sufficient condition for the finiteness conjecture on $M=G/K$ to hold: any $G$-invariant metric is non-degenerate.  On the other hand, if $g$ is {\it unstable} (i.e., $\scalar''_g(T,T)>0$ for some $T\in{\tca\tca_g^G}$), then $g$ is also unstable relative to the $\nu$-entropy functional introduced by Perelman, it is dynamically unstable and it does not realize the Yamabe invariant of M.

We study in this paper the stability type of Einstein metrics on flag manifolds $M=G/T$, where $G$ is a compact simple Lie group and $T$ a maximal torus of $G$.  Note that $\Aut(G/K)$ is therefore always finite.  The finiteness conjecture is still open for most flag manifolds, the simplest cases for which we do not know whether there are infinitely many Einstein metrics or not are $M^{12}=\SO(6)/T^3$ and $M^{20}=\SU(5)/T^4$, which have $\dim{\mca^G_1}=5$ and $9$, respectively.  

As an application of the method given in \cite{LauretWill}, we first compute the stability type, including the coindices, of all invariant Einstein metrics on 
$$
G_2/T^2, \qquad \SO(5)/T^2, \qquad \Spe(3)/T^3, \qquad \SU(3)/T^2, \qquad \SU(4)/T^3,
$$
and of all known ones on $\SO(6)/T^3$ (see Tables 1-6).  

Secondly, we focus in Section 7 on the homogeneous space 
$$
F(5):=\SU(5)/T^4.
$$  
We explicitly give four new Einstein metrics on $F(5)$, which were numerically predicted in \cite{linoGrama}.  Together with the standard metric, the K\"ahler metric and the Arvanitoyeorgos metric, this gives a total of seven Einstein metrics on $F(5)$.  We obtain that they are all non-degenerate and unstable, but surprisingly, the approximate (volume normalized) scalar curvature values of these seven Einstein metrics are all in a range of less than $0.05$, so this invariant becomes useless to differentiate them up to homothecy.  However, we show that the stability data, including the coindex and some information on the spectrum of the Hessian $\scalar''_g$, is enough to obtain that these seven Einstein metrics on $F(5)$ are pairwise non-homothetic.

\section{Preliminaries}

Let $M$ be a compact homogeneous connected differentiable manifold, we fix an almost-effective transitive action of 
a compact Lie group $G$ on $M$. This action determines a presentation $M=G/K$, where $K$ is the isotropy subgroup at some point $o\in M$. 
Let $\mca^G$ be the space of all $G$-invariant metrics on $M$ and $\mca^G_1$ the space of all $G$-invariant unit volume metrics. 

We say that a Riemannian manifold $(M,g)$ 
is Einstein if $Rc(g)=\rho g$ for some $\rho\in\RR$.  We will focus on $G$-invariant metrics, an alternative characterization of Einstein $G$-invariant metrics is that they are precisely the critical points of the scalar curvature functional
$$
Sc:\mathcal{M}_1^G\to \RR.
$$
The signature of the Hessian $Sc''$ determines the stability type of the metric.  The tangent space of $\mca_1^G$ at $g$ is the $g$-traceless subspace of $S^2(M)^G$, where
$S^2(M)^G$ is the space of $G$-invariant symmetric 2-tensors.

It is well-known that 
$$
T_g\mca_1^G=T_g\Aut(G/K)\cdot g\oplus\tca\tca_g^G,
$$
where $T_g\Aut(G/K)\cdot g$ is the space of trivial variations of $g$ and $\tca\tca_g^G$ is the space of $G$-invariant traceless transversal tensors, called TT-tensors. We will be interested on the coindex of $g$, which is the dimension of the maximal subspace of $\tca\tca_g^G$ where $Sc_g''$ is positive definite. 
\begin{definition}
    An Einstein metric $g$ is called
    \begin{itemize}
        \item $G$-stable if $\left. Sc''(g)\right|_{\tca\tca_g^G}<0$.
        \item $G$-unstable if $\left. Sc''(g)\right|_{\tca\tca_g^G}(T,T)>0$ for some $T\in\tca\tca_g^G$. 
        \item $G$-non-degenerate if $\left. Sc''(g)\right|_{\tca\tca_g^G}$ is non-degenerate. 
    \end{itemize}
\end{definition}

We know that if $M=G/K$ is as above, then there is at least one reductive decomposition $\ggo=\mathfrak{k}\oplus\pg$.
We can identify $\pg$ with $T_oM$, and $\mathfrak{k}$ is the Lie algebra of $K$. We extend $\langle\cdot,\cdot\rangle:=g_o$ to $\glg(\pg)$ and $\Lambda^2\pg^*\otimes \pg$ as follows:
$$
\langle A,B\rangle:=\operatorname{tr}AB^t\quad\text{and}\quad \langle \lambda,\lambda\rangle:=\sum|\lambda(X_i,X_j)|^2,
$$
where $A,B\in\glg(\pg)$, $\lambda\in \Lambda^2\pg^*\otimes \pg$ and $\{X_i\}$ is any orthonormal basis of $\pg$.

\subsection{Stability in terms of the Lichnerowicz Laplacian.}
Let $(g,\rho)$ be an Einstein metric on $M$. A classic result is the relation between the Lichnerowicz Laplacian and the Hessian, given by
$$
\left.Sc''_g\right|_{\tca\tca_g}=\tfrac{1}{2}\langle (2\rho\operatorname{id}-\Delta_L)\cdot,\cdot \rangle,
$$
where $\Delta_L$ is the Lichnerowicz Laplacian of $g$. 

Consider any reductive decomposition $\ggo=\mathfrak{k}\oplus\pg$ and the usual identifications $\pg\equiv T_oM$ and $\mathcal{S}^2(M)^G \equiv \operatorname{sym}(\mathfrak{p})^K$, where $\operatorname{sym}(\mathfrak{p})^K$ is the vector space of all $Ad(K)$-invariant symmetric linear transformations on the
$n$-dimensional vector space $\pg$. We denote by $\mu$ the Lie bracket of $\ggo$ and $\mu_\pg:\pg\times\pg\to\pg$ is given by $\left.\operatorname{pr}_\pg \circ \mu\right|_{\pg\times\pg}$, where $\operatorname{pr}_\pg$ is the projection on $\pg$ relative to $\ggo=\kg\oplus\pg$. Let $M_{\mu_\pg}$ be the operator defined by
$$
\langle M_{\mu_\pg}, A\rangle=\langle \theta (A) \mu_{\pg},\mu_\pg\rangle,
$$
here $\theta$ is given by
$$
\theta (A) \lambda:= A\lambda(\cdot,\cdot)-\lambda(A\cdot,\cdot)-\lambda(\cdot,A\cdot).
$$

Consider the operator $\Le_\pg=\Le_\pg(g):\operatorname{sym}\pg\to\operatorname{sym}\pg$ defined by 
$$
\langle \Le_\pg A,B\rangle=\tfrac{1}{2}\langle\theta(A)\mu_\pg,\theta(B)\mu_\pg\rangle+2\operatorname{tr}M_{\mu_\pg} AB \qquad\forall A,B\in \operatorname{sym}(\pg),
$$
where $\operatorname{sym}(\pg)$ denotes the space of symmetric linear transformations on $\pg$.
In \cite[Corollary 4.7]{jLauret}, the author proved that under identifications, $\Delta_LT=\langle \Le _\pg A\cdot,\cdot\rangle,$ where $T\in\tca\tca_g^G$ and $T=\langle A\cdot,\cdot\rangle$ for some $A\in\operatorname{sym}(\pg)^K$.

According to this result, the stability of $g$ can be expressed in terms of the minimum and maximum eigenvalues of $\left.\Le_\pg\right|_{\tca\tca_g^G}$, which will be denoted by $\lambda_\pg$ and $\lambda_\pg^{\max{}}$, respectively, as follows:

\begin{itemize}
    \item $G$-stable if and only if $2\rho<\lambda_\pg$.
    \item $G$-unstable if and only if $\lambda_\pg<2\rho$.
    \item $G$-non-degenerate if and only if $2\rho\not\in\operatorname{Spec}\left(\left.\Le_\pg \right|_{\tca\tca_g^G}\right)$.
    \item A local minimum of $\left.Sc\right|_{\mca_1^G}$ if and only if $\lambda_\pg^{\operatorname{max}}<2\rho$.    
\end{itemize}
Also we can write the coindex of $g$ in terms of the spectrum of $\Le_\pg$, the coindex of $g$ is equal to $(\sum_{\lambda<2\rho}\dim\operatorname{Eigenspace}(\lambda))-1$,
where the sum is taken over the spectrum of $[\Le_\pg]$. 

Therefore, the study of the stability of Einstein metrics is reduced to the study of the spectrum of $\Le_\pg(g)$.

\subsection{Structural constants.}

Let $Q$ be a bi-invariant inner product on $\ggo$, we fix a $Q$-orthogonal reductive decomposition 
$$
\ggo=\mathfrak{k}\oplus\pg_1\oplus\cdots\oplus\pg_r,
$$
in $\operatorname{Ad}(K)$-invariant subspaces. Let $\{e_\alpha\}$ be an $Q$-orthogonal basis of $\pg$, adapted to the above decomposition.
That is, $e_\alpha\in \pg_i$ for some $i$, and if $e_\alpha\in \pg_i,e_\beta\in\pg_j$ and $i<j$ then $\alpha<\beta$. We define the structural coefficients $[ijk]$ due to \cite[\S 1]{wangZillerSoloCoeficientesEstructurales}, as follows:

$$
[ijk]=\sum_{\alpha,\,\beta,\,\gamma}Q([e_\alpha^i,e_\beta^j],e_\gamma^k)^2,
$$
where $\{e^s_\alpha\}$ is a $Q$-orthonormal basis of $\pg_s$. As we can see, $[ijk]$ is invariant under any permutation of the indices. 
Also, $[ijk]=0$ if and only if $[\pg_i,\pg_j]$ is orthogonal to $\pg_k$.

For each $g\in\mathcal{M}$, there exists at least one decomposition 
\begin{equation}\label{dec}
\pg=\pg_1\oplus\cdots\oplus\pg_r,    
\end{equation}
which is also $g$-orthogonal and, 
$$
g=x_1\left.Q\right|_{\pg_1}+\cdots+x_r\left.Q\right|_{\pg_r},
$$
note that necessarily $x_i>0$. The metric $g$ will be denoted by $(x_1,...,x_r)$.

The structural constants allow us to express the Einstein equations in terms of the $x_i$'s.
If $G$ is a simple Lie group, $Q=-B_\ggo$ the negative of the Killing's form and the decomposition \eqref{dec} is multiplicity free, we can write the Ricci eigenvalues as follows:
\begin{equation}\label{eq:Einstein Eq}
\rho_k=\tfrac{1}{2x_k}+\tfrac{1}{4d_k}\sum_{i,j}\tfrac{x_k^2-x_i^2-x_j^2}{x_ix_jx_k}[ijk],\qquad k=1,\ldots,r. 
\end{equation}
Another important result is \cite[Theorem 3.1.]{LauretWill}, which states that $\Le_{\pg}$ is also in terms of the $x_i$'s, we consider the orthonormal basis of $\operatorname*{sym}(\pg)^K$ given by
$$
\left\{\tfrac{1}{\sqrt{d_1}} I_1, \ldots, \tfrac{1}{\sqrt{d_r}} I_r\right\},
$$
where $I_k$ is the block map given by $[0,..,0,I_{\pg_k},0,...,0]$ of $\pg$. We have the following result:
\begin{equation}\label{eq:L}
    \begin{aligned}
    \relax [\Le_\pg]_{kk}&=\tfrac{1}{d_k}\sum_{i,j\neq k}\tfrac{x_k}{x_ix_j}[ijk]+\tfrac{1}{d_k}\sum_{i\neq k}\tfrac{x_i}{x_k^2}[ikk],\\
    [\Le_\pg]_{km}&=\tfrac{1}{\sqrt{d_k}\sqrt{d_m}}\sum_{i}\tfrac{x_i^2-x_k^2-x_m^2}{x_ix_kx_m}[ikm].
    \end{aligned}
\end{equation}
Finally we can write the volume normalized scalar curvature in the same terms as follows:
$$
Sc_N(g):=(\operatorname{det} _Q(g))^{\tfrac{1}{n}}\operatorname{Sc}(g)=(x_1^{d_1}\cdots x_r^{d_r})^{\tfrac{1}{n}}\left(\tfrac{1}{2}\sum_k\tfrac{d_k}{x_k}-\tfrac{1}{4}\sum_{i,j,k}\tfrac{x_k}{x_ix_j}[ijk]\right).
$$
\begin{remark}
    The value of $Sc_N$ is an homothetic invariant.
\end{remark}
To end this subsection, we will give another motivation for the study of $G$-stability in this case.
\begin{definition}
    A $G$-invariant Einstein metric $g$ is called $G$-rigid if there exists an open
neighborhood $U$ of $g$ in $\mca^G$ such that any Einstein $g_0\in U$ belongs to $\operatorname{Aut}(G/K)\cdot g$ up to
scaling.
\end{definition}

Let $\mathcal{E}^G_1$ be the space of $G$-invariant Einstein metrics,  so we got the following result due to \cite[Theorem 1.6]{BWZSoloTeoremaECompacto}: 
\begin{theorem}\label{thm: BWZ}
    Let $G$ be a compact Lie group and $M = G/K$ be a connected homogeneous space with finite fundamental group. Then each connected component of $\mathcal{E}^G_1$ is compact and the set of possible Einstein constants of metrics among $\mathcal{E}^G_1$ is finite.
\end{theorem}
Then the space $\mathcal{E}^G_1/\operatorname*{Aut}(G/K)$ is also compact, and hence $\mathcal{E}^G_1/\operatorname*{Aut}(G/K)$ is finite if and only if every $g \in \mathcal{M}^G$ is $G$-rigid. 
A metric being G-non-degenerate implies G-rigidity. Therefore, if every metric $g\in \mca^G$ is non-degenerate then $\mathcal{E}^G_1/\operatorname*{Aut}(G/K)$ is finite. 

In the following sections, we will describe the $G$-stability of some full-flag manifolds.

\section{Einstein Metrics on \texorpdfstring{$G_2/T^2$}{G2/T2}}
\section{Einstein Metrics on \texorpdfstring{$SO(5)/T^2$}{SO(5)/T2}}
In this section, we will investigate Einstein metrics on $SO(5)/T^2$. In this particular scenario, $\dim \mathcal{M}_1^G=3$ and all metrics are precisely documented. The isotropy representation is given by 
$$
\pg=\pg_1\oplus\pg_2\oplus\pg_3\oplus\pg_4,
$$
and
$$
[1\,3\,4]= [2\,3\,4]=\frac{1}{3}, \qquad \text{zero otherwise.}
$$
In \cite{Sakane} the author shows that, up to scaling, there are exactly 6 solutions to the Einstein equation,
$$
\begin{aligned}
g_{K_1}&=(2,4,1,3), \quad g_{K_2}=(4,2,1,3), \quad g_{K_3}=(2,4,3,1), \quad g_{K_4}=(4,2,3,1),\quad \rho(g_{K_i})=\tfrac{1}{6},\\
g_1&=(\tfrac{24+4\sqrt{6}}{15},\tfrac{24+4\sqrt{6}}{15},1,\tfrac{7+2\sqrt{6}}{5}),\qquad \rho(g_1)=\tfrac{6-\sqrt{6}}{18}, \\
g_2&=(\tfrac{24-4\sqrt{6}}{15},\tfrac{24-4\sqrt{6}}{15},1,\tfrac{7-2\sqrt{6}}{5}),\qquad\rho(g_1)=\tfrac{6+\sqrt{6}}{18},
\end{aligned}
$$
we got the following result.
\begin{theorem}
The full flag manifold $SO(5)/T$ admits exactly two $G$-invariant Einstein metrics (up to homothecy). These metrics are the Kähler metric $g_K=g_{K_1}$ and $g_1$, their $G$-stablity is given in Table \ref{tab:metrics SO5/T}.
\end{theorem}
\begin{proof}
    Since $ g_K $ and $ g_1 $ have different coindices, they are not homothetic. It remains to show that $ g_1 $ and $ g_2 $ are homothetic. To see this, observe that
$\rho(g_1) g_1 = \rho(g_2) g_2.$
    
\end{proof}

\begin{table}
    \centering
\begin{tabular}{| c | c | c | c | c | c | c | c |}
\hline
Name  & Propieties & $2\rho$ & $Sc_N$ & coindex & $\lambda_\pg$ &$\lambda_{\pg}^{max}$\\ \hline
$g_K$  & Non-deg. and unstable. & {\tiny$\tfrac{1}{6}\approx 1.6667$} & $2.9511$ & 1 & 0 & $1.0654$\\ \hline
$g_1$  & Non-deg. and unstable. & $0.4694$ & $2.9420$
 & 2 & $0.3755$ & $2.3470$\\ \hline
\end{tabular}
    \caption{$G$-stability of metrics on $SO(5)/T^2$.}
    \label{tab:metrics SO5/T}
\end{table}

\section{Einstein Metrics on \texorpdfstring{$SO(6)/T^3$}{SO(6)/T3}}
In this section, we will investigate Einstein metrics on $SO(6)/T$. Following \cite[Section 3.4]{Sakane} the root representation is given by 

$$
\begin{aligned}
\pg&=\pg_{\varepsilon_1-\varepsilon_2}\oplus\pg_{\varepsilon_1-\varepsilon_3}\oplus\pg_{\varepsilon_2-\varepsilon_3}\oplus\pg_{\varepsilon_1+\varepsilon_2}\oplus\pg_{\varepsilon_1+\varepsilon_3}\oplus\pg_{\varepsilon_2+\varepsilon_3}\\
&=\pg_{1}\oplus\pg_{2}\oplus\pg_{3}\oplus\pg_{4}\oplus\pg_{5}\oplus\pg_{6},
\end{aligned}
$$
therefore $\dim\mathcal{M}_1^G=5$. Furthermore,
$$
[1\,2\,3]= [1\,5\,6]=[2\,4\,6]=[3\,4\,5]=\frac{1}{4}, 
$$
and otherwise zero, there are at least two Einstein metrics, given by
$$
g_{S}=(1,1,1,1,1,1) \qquad g_1=(5,5,5,3,3,3).
$$

\begin{theorem}
The full flag manifold $SO(6)/T^3$ admits at least two $SO(6)$-invariant Einstein metrics up to homothecy and their $G$-stablity is given in  Table \ref{tab:metrics SO6/T}.
\end{theorem}
\begin{proof}
    We only have to prove that these metrics are not homothetic, observe that their coindex are different.
\end{proof}
\begin{remark}
    In \cite[Section 5]{LauretWillLocalMaxima} the authors proved that the standar metric  on $SU(2n)/T$ is actually a saddle point.
\end{remark}
\begin{table}
    \centering
\begin{tabular}{| c | c | c | c | c | c | c | c |}
\hline
Name  & Propieties & $2\rho$ & $Sc_N$ & coindex & $\lambda_\pg$ &$\lambda_{\pg}^{max}$\\ \hline
$g_S$  & Deg. and unstable (saddle) & {\tiny$\tfrac{3}{4}$} & {\tiny$\tfrac{9}{4}=2.25$} & {\tiny 3} & {\tiny$\tfrac{1}{2}$} & {\tiny$\tfrac{3}{4}$}\\ \hline
$g_1$  & Non-deg. and unstable. & {\tiny$\tfrac{7}{36}$} & {\tiny$\tfrac{11}{20}\sqrt{15}\approx 2.1301$} & {\tiny 2} & {\tiny$\tfrac{53}{360}-\tfrac{1}{360}\sqrt{1201}$} & {\tiny$\tfrac{5}{18}$}\\ \hline
\end{tabular}
    \caption{$G$-stability of metrics on $SO(6)/T^3$.}
    \label{tab:metrics SO6/T}
\end{table}

\section{Einstein Metrics on \texorpdfstring{$Sp(3)/T^3$}{Sp(3)/T3}}
In this section, we will explore the Einstein metrics on $Sp(3)/T^3$. Specifically, all metrics pertinent to this scenario have been found in \cite{GAO-WANG-Sp3}, we will show their $G$-stability. The isotropy representation is given by
$$
\mg=\mg_1\oplus\mg_2\oplus\mg_3\oplus\mg_4\oplus\mg_5\oplus\mg_6\oplus\mg_7\oplus\mg_8\oplus\mg_9,
$$
where $\dim \mg_i=2$. Threfore $\dim\mathcal{M}_1^G=8$ and
$$
\begin{aligned}
\relax [1\,2\,4]&=[2\,4\,5]=[1\,6\,7]=[3\,5\,8]=[3\,8\,9]=[6\,7\,9]=\frac{1}{4}, \\
[2\,3\,6]&=[3\, 4\, 7]=[4\, 6\, 8]=[2\, 7\, 8]=\frac{1}{8}.
\end{aligned}
$$

\begin{theorem} \cite[Theorem 2]{GAO-WANG-Sp3} 
The full flag manifold $Sp(3)/T$ admits exactly four $Sp(3)$-invariant Einstein metrics. These metrics are specified as follows:
$$
\begin{aligned}
g_K&=(1,\tfrac{1}{4},\tfrac{1}{2},\tfrac{3}{4},\tfrac{1}{2},\tfrac{1}{4},\tfrac{5}{4},1,\tfrac{3}{2}),\\ 
g_1&\approx(1,0.4311,1.0381, 1.0381, 1, 1.0381, 0.4311, 0.4312, 1),\\
g_2&\approx(1, 0.3430, 0.9326, 0.9326, 0.8101, 1.0567, 0.5477, 0.3430,1),\\
g_3&\approx (1, 0.3149, 0.8524, 0.9092, 0.7740, 0.9896, 0.5185, 0.3298, 0.8708),
\end{aligned}
$$ 
where the first one is Kähler, and the last three are given by approximate values.
\end{theorem}

\begin{theorem}
The metrics $g_K, g_1$, $g_2$ and $g_3$  are pairwise non-homothetic, and their $G$-stability is presented in Table \ref{tab:metrics Sp3/T}.
\end{theorem}
\begin{proof}
Observe that $g_1$ has coindex 3, while $g_K$ and  $g_2$ have coindex 2, implying that $g_1$ is non-homothetic to both $g_K$ and $g_2$. Additionally,  $\lambda_\pg(g_K) = 0$ whereas $\lambda_\pg(g_2) \neq 0$, making $g_K$ and $g_2$ non-homothetic. Finally, we have
$$
\frac{\lambda_\pg^{\max}(g_1)}{\lambda_\pg(g_1)} = 6.050 \neq 8.9073 = \frac{\lambda_\pg^{\max}(g_3)}{\lambda_\pg(g_3)},
$$
confirming that all metrics are non-homothetic.
\end{proof}

\begin{table}
    \centering
\begin{tabular}{| c | c | c | c | c | c | c | c |}
\hline
Name  & Propieties & $2\rho$ & $Sc_N$ & coindex & $\lambda_\pg$ &$\lambda_{\pg}^{max}$\\ \hline
$g_K$  & Non-deg. and unstable. & {\tiny$\tfrac{1}{2}$} & 5.8885 & 2 & 0 & 1.028\\ \hline
$g_1$  & Non-deg. and unstable. & 0.8528 & 5.8711 & 3 & 0.3790 & 2.2931\\ \hline
$g_2$  & Non-deg. and unstable. & 0.9149 & 5.8759 & 2 & 0.3552 & 2.8374\\ \hline
$g_3$  & Non-deg. and unstable. & 0.9719 & 5.8759 & 3 & 0.3420 & 3.0463\\ \hline
\end{tabular}
    \caption{$G$-stability of metrics on $Sp(3)/T^3$.}
    \label{tab:metrics Sp3/T}
\end{table}
\begin{remark}
    The values with a decimal point are approximate values.
\end{remark}
\section{Einstein Metrics on \texorpdfstring{$F(n)$}{F(n)}}
Our aim in this section is to study the case 
$$
F(n):=SU(n)/T^{n-1}.
$$ 

The reductive decomposition is given by:
\begin{equation}{\label{reductive decomposition}}
\mathfrak{su}(n)=\mathfrak{t}\oplus\pg_{12}\cdots\oplus\pg_{ij}\oplus\cdots\pg_{(n-1)n}.
\end{equation}
 Each of these $p_{ij}$ is irreducible and they are all pairwise non-equivalent, meaning that $F(n)$ is multiplicity free. 
 Another important property that these $\pg_{ij}$ satisfy is that $[\pg_{ij},\pg_{kl}]_{\pg}=0$ if $\{i,j\}$ and $\{k,l\}$ are disjoint or equal, and $[\pg_{ij},\pg_{ik}]_{\pg}$ is nonzero and contained in $\pg_{jk}$. 
 So the only nonzero coefficients are $[(ij)(jk)(ik)]$, and it is wellknown that its value is $\tfrac{1}{n}$, see \cite[Section 3.1.]{Sakane} for more details.

\begin{proposition}\label{prop: eq for F(n)}
Given a metric $g=(x_{12},..,x_{(n-1)n})$ on $F(n)$, the formulas reduce to:
$$
\begin{aligned}
    2\rho_{ij}&=\tfrac{1}{x_{ij}}+\tfrac{1}{2n}\sum_{k\neq i,j}\tfrac{x_{ij}^2-x_{ik}^2-x_{jk}^2}{x_{ij}x_{ik}x_{jk}},\\
    [\Le_\pg]_{(ij)(ij)}&=\tfrac{1}{n}\sum_{k\neq i,j}\tfrac{x_{ij}}{x_{ik}x_{jk}},
    \quad\quad [\Le_\pg]_{(ij)(jk)}=\tfrac{1}{2n}\tfrac{x_{ik}^2-x_{ij}^2-x_{jk}^2}{x_{ij}x_{ik}x_{jk}}.
\end{aligned}
$$
The remaining entries of $[\Le_\pg]$ are zero. 
\end{proposition}
\begin{proof}
    Let us start with $\rho_{ij}$. Remembering that $d_{ij}=2$ and that the only nonzero structural constants are $[(ij)(jk)(ik)]$, we conclude that in the sum over $((km),(rt))$, the only nonzero elements in the sum are when the pairs $((km),(rt))=((j\alpha),(i\alpha))$ and $((km),(rt))=((i\alpha),(j\alpha))$. 
    With this clarified, we can start the computations, using \eqref{eq:Einstein Eq} we get:
    $$
\begin{aligned}
2\rho_{ij}-\tfrac{1}{x_{ij}}&=\tfrac{1}{2d_k}\sum_{(km),(rt)}\tfrac{x_{ij}^2-x_{km}^2-x_{rt}^2}{x_{ij}x_{km}x_{rt}}[(ij)(km)(rt)]\\
&=\tfrac{1}{4}\left[\sum_{\alpha}\tfrac{x_{ij}^2-x_{j\alpha}^2-x_{i\alpha}^2}{x_{ij}x_{j\alpha}x_{i\alpha}}[(ij)(j\alpha)(i\alpha)]+\sum_\alpha\tfrac{x_{ij}^2-x_{i\alpha}^2-x_{j\alpha}^2}{x_{ij}x_{i\alpha}x_{j\alpha}}[(ij)(i\alpha)(j\alpha)]\right]\\
&=\tfrac{1}{2}\sum_{\alpha}\tfrac{x_{ij}^2-x_{j\alpha}^2-x_{i\alpha}^2}{x_{ij}x_{j\alpha}x_{i\alpha}}[(ij)(j\alpha)(i\alpha)]\\
&=\tfrac{1}{2n}\sum_{\alpha}\tfrac{x_{ij}^2-x_{j\alpha}^2-x_{i\alpha}^2}{x_{ij}x_{j\alpha}x_{i\alpha}}.
\end{aligned}
$$
 We observe that $[(km)(ij)(ij)]=0$, now we can compute $[\Le_\pg]_{(ij)(ij)}$:
    
    $$ 
    \begin{aligned}
   \relax [\Le_\pg]_{(ij)(ij)}&=\tfrac{1}{d_{ij}}\sum_{(km),(rt)\neq(ij)}\tfrac{x_{ij}}{x_{km}x_{rt}}[(ij)(km)(rt)]+\tfrac{1}{d_{ij}}\sum_{(km)\neq (ij)}\tfrac{x_{km}}{x_{ij}^2}[(km)(ij)(ij)]\\
    &=\tfrac{1}{2}\sum_{(km),(rt)\neq(ij)}\tfrac{x_{ij}}{x_{km}x_{rt}}[(ij)(km)(rt)]\\
    &=\tfrac{1}{2}\left[\sum_\alpha\tfrac{x_{ij}}{x_{j\alpha}x_{i\alpha}}[(ij)(j\alpha)(i\alpha)]+
    \sum_\alpha\tfrac{x_{ij}}{x_{i\alpha}x_{j\alpha}}[(ij)(i\alpha)(j\alpha)]\right]\\
    &=\tfrac{1}{n}\sum_\alpha\tfrac{x_{ij}}{x_{j\alpha}x_{i\alpha}}.
    \end{aligned}
    $$
Let us continue with $[\Le_\pg]_{(ij)(jk)}$:
    $$
\begin{aligned}
\relax [\Le_\pg]_{(ij)(jk)}&=\tfrac{1}{\sqrt{d_{ij}}\sqrt{d_{jk}}}\sum_{(rt)}\tfrac{x_{rt}^2-x_{ij}^2-x_{jk}^2}{x_{ij}x_{jk}x_{rt}}[(ij)(jk)(rt)]\\
&=\tfrac{1}{\sqrt{2}\sqrt{2}}\tfrac{x_{ik}^2-x_{ij}^2-x_{jk}^2}{x_{ij}x_{jk}x_{ik}}[(ij)(jk)(ik)]\\
&=\tfrac{1}{2n}\tfrac{x_{ik}^2-x_{ij}^2-x_{jk}^2}{x_{ij}x_{jk}x_{ik}}.
\end{aligned}
$$
Finally, the remaining terms of the matrix $[\Le_\pg]_{(ij)(km)}$ are zero because the coefficient $[(ij)(km)(rt)]$ appears in all elements of the summation, which is zero since $\{i,j\}\cap\{k,m\}$ is empty.
\end{proof}

\begin{remark}
    It follows from the previous proposition that
    $$
    2\rho_{ij}=\tfrac{1}{x_{ij}}+\sum_{k\neq i,j}[\Le_\pg]_{(ik)(jk)}.
    $$
\end{remark}

There are some families of Eintein metrics on $F(n)$ known in the literature, for example the following:
\begin{itemize}\label{classic-etrics}
    \item For $n\geq 3$, let $x_{ij}=1$ for all $i,j$. This metric is called the \textit{standard metric}.
    \item For $n\geq 3$, let $x_{1k}=n-1$, for $k=2,..,n$ and $x_{ij}=n+1$, $1<i<j$. This metric is called the \textit{Arvanitoyeorgos metric} \cite[Theorem 7]{Arvanitoyeorgos}.
    \item For $n=2m$, $m\geq 3$, let $x_{ij}=m+2$ for $1\leq i<j\leq m$ or $m<i<j\leq n$ and $x_{ij}=3m-2$ otherwise. This metric is called the \textit{Senda metric} \cite{senda1997}. 
\end{itemize}
For the above metrics, we have proved something stronger:
\begin{proposition}\label{prop:unicidad particion}
    \begin{enumerate}        
        The following statements are true:
        \item On $F(n)$, the only Einstein metrics such that $x_{1i}=A$ for $i=2,...,n$ and $x_{ij}=B$ otherwise are the standard metric and the Arvanitoyeorgos metric.
        \item On $F(2m)$, the only Einstein metrics such that $x_{ij}=A$ for $1\leq i<j\leq m$, $x_{ij}=C$ for $m<i<j\leq n$ and $x_{ij}=B$ otherwise are the standard metric and the Senda metric.
    \end{enumerate}
\end{proposition}
\begin{proof}
    We will only prove part (ii), since part (i) is completely analogous. To prove this we will compute $\rho_{ij}$ for the distinct cases of $i,j$.
    First case $1\leq i<j\leq m$, then a straightforward computation using Proposition \ref{eq:Einstein Eq} shows that if $i,j\leq m$ then
    $$
    \rho_{ij}=\tfrac{3m+2}{4m}\tfrac{1}{A}+\tfrac{1}{4}\tfrac{A^2-2B^2}{AB^2},\\   
    $$
    if $m<i,j$, then
    $$
    \rho_{ij}=\tfrac{3m+2}{4m}\tfrac{1}{C}+\tfrac{1}{4}\tfrac{C^2-2B^2}{CB^2},\\
    $$
    and otherwise
    $$
    \rho_{ij}=\tfrac{1}{B}-\tfrac{m-1}{4m}\left(\tfrac{A}{B^2}\right)-\tfrac{m-1}{4m}\tfrac{C}{B^2}.\\
    $$
    The metric is Einstein if the above expressions are equal, which is equivalent to the following system of equations:
    $$
    \left\lbrace\begin{aligned}
        &(A-C)(mCA-(m+2)B^2)=0\\ 
        &(m+2)B^2+(2m-1)A^2-4mAB+(m-1)AC=0 .
        \end{aligned}\right.        
    $$    
The case $A-C=0$ gives us the standard metric and the Senda metric as solutions. If $A-C\neq 0$, then necessarily $mCA-(m+2)B^2=0$, which does not have real solutions.
\end{proof}

\begin{remark}
In \cite{Negreiros}, the authors claim three family of solutions to the Einstein equations, which we unfortunately proved that they are not actually Einstein. Such metrics are the following
\begin{itemize}
    \item If $n=2m$ and $m\geq 6$ with $x_{ij} = m + 5$ for $1\leq i<j\leq m$ or $m<i<j\leq n$ and  $x_{ij} = 3m - 5$ otherwise. 
    \item If $n=2m+1$, $m\geq 6$ with $x_{ij} = A$ for $1\leq i<j\leq m+1$ or $m+1<i<j\leq n$ and $x_{ij} = B$ otherwise.
    \item On $F(5)$ take $(1,1,2,2,1,2,2,1,1)$. 
\end{itemize}
\end{remark}
 
\subsection{Einstein Metrics on \texorpdfstring{$F(3)$}{F(3)} and \texorpdfstring{$F(4)$}{F(4)}}
We now focus on the particular cases $F(3)$ and $F(4)$, with $\dim \mathcal{M}_1^G=2$ and $\dim \mathcal{M}_1^G=5$, respectively. We will compute the $G$-stability types of the Eintein metrics found in the literature. 
\subsection*{Case  \texorpdfstring{$F(3)$}{F(3)}}
This case has been studied in \cite[Section 4.1.]{LauretWill}, as a particular case of generalized Wallach spaces, we could observe $\dim\mathcal{M}_1^G=2$. Up to scaling and isometries, there are only two Einstein metrics on $F(3)$, 
$$
g_s=(1,1,1) \qquad \text{and} \qquad g_A=(1,1,2).
$$
The standard metric $g_s$ with the Einstein constant $\rho=\tfrac{5}{12}$ is a local minimum, in particular is non-degenerate and unstable, with coindex 2. On the other hand, the Arvanitoyeorgos metric $g_A$, with $\rho=\tfrac{1}{3}$ is non-degenerate, unstable with coindex 1, and in this particular scenario is a K\"aler metric.

\begin{table}
    \centering
\begin{tabular}{| c | c | c | c | c | c | c | c |}
\hline
Name  & Propieties & $2\rho$ & $Sc_N$ & coindex & $\lambda_\pg$ &$\lambda_{\pg}^{max}$\\ \hline
$g_S$  & local min. & \small{$\tfrac{5}{6}\approx 0.8333$} & \small{$\tfrac{5}{2}$} & 2 & \small{$\tfrac{1}{2}$} & \small{$\tfrac{1}{2}$}\\ \hline
$g_A$  & Non-deg. and unstable. & \small{$\tfrac{2}{3}$} & 2.5198 & 1 & 0 & 1\\ \hline
\end{tabular}
    \caption{$G$-stability of metrics on $SU(3)/T$.}
    \label{tab:metrics SU3/T}
\end{table}

\subsection*{Case  \texorpdfstring{$F(4)$}{F(4)}}
This particular scenario has been studied in \cite{Wang_Li_2014_F4},  the authors claim that there are only four Einstein metrics up to homothety. We have found their expressions in our study, the metrics are the standard metric $g_s$, the Arvanitoyeorgos metric $g_A$, the K\"ahler metric $g_K$ and $g_1$. Where
$$
\begin{aligned}
&g_s=(1,1,1,1,1,1),\qquad g_A=(3,3,3,6,6,6), \\
&g_K=(1,2,3,1,2,1),\qquad g_1\approx(1,1.9436,1.9436,1.1867,1.1867,1.2815).
\end{aligned}
$$
The following result demonstrates that their $G$-stability type distinguishes them.
\begin{theorem}
The full flag manifold $SU(4)/T$ admits exactly four $G$-invariant Einstein metrics (up to homothecy). These metrics are the Kähler metric $g_K$, the standar metric $g_s$, the Arvanitoyeorgos metric $g_A,$ and $g_1$. Their stability type is given in \ref{tab:metrics SU4/T}.
\end{theorem}
\begin{proof}
    Observe that by coindex, we only need to prove that $g_A$ is not homothetic to $g_K$ and $g_S$ is not homothetic to $g_1$. For the first, $\lambda_\pg(g_K)$ is zero and $\lambda_\pg(g_A)\neq 0$. For the second, $g_S$ is $G$-degenerate and $g_1$ $G$-nondegenerate.
\end{proof}

\begin{table}
    \centering
\begin{tabular}{| c | c | c | c | c | c | c | c |}
\hline
Name  & Propieties & $2\rho$ & $Sc_N$ & coindex & $\lambda_\pg$ &$\lambda_{\pg}^{max}$\\ \hline
$g_S$  & deg. and unstable. & \small{$\tfrac{3}{4}$} & 4.5 & 3 & \small{$\tfrac{1}{2}$} & \small{$\tfrac{3}{4}$}\\ \hline
$g_A$  & Non-deg. and unstable. & \small{$\tfrac{7}{36}$} & 4.5184 & 2 & 0.0509 & 0.2778\\ \hline
$g_K$  & Non-deg. and unstable. & \small{$\tfrac{1}{2}$} & 4.5392 & 2 & 0 & 1.1371\\ \hline
$g_1$  & Non-deg. and unstable. & 0.5462 & 4.5136 & 3 & 0.0875 & 0.8731\\ \hline
\end{tabular}
    \caption{$G$-stability of metrics on $SU(4)/T$.}
    \label{tab:metrics SU4/T}
\end{table}

Given that the standard metric is $G$-degenerate (see \cite{jLauret}), it can either be a local minimum or a saddle point. We have shown that the standard metric is, in fact, a saddle point.
\begin{proposition}
The standard metric on $F(4)$ is a saddle point.    
\end{proposition}
\begin{proof}
First, consider a metric $(x_{12},x_{13},x_{14},x_{23},x_{24},x_{34})$, and take $x_{12}=x_{13}=X,x_{14}=x_{23}=Y$, $x_{24}=x_{34}=Z$, and $Z=\tfrac{1}{XY}$, then
{\small
$$
Sc_N(X,Y)=-\tfrac{1}{8}\,{\tfrac {2\,{Y}^{2}{X}^{4}-12\,{X}^{2}Y+2\,{Y}^{4}{X}^{2}+{Y}^{5}{X}^{4}+{Y}^{3}+2-16\,{X}^{3}{Y}^{3}-16\,{Y}^{2}X}{{Y}^{2}{X}^{2}}}.
$$
}
If we compute the hessian matrix of $Sc_N(X,Y)$ on $(1,1)$, then we see that the null space is generated by $(-\tfrac{1}{2},1)$, this motivates us to set $X=-\tfrac{1}{2}Y+\tfrac{3}{2}$. Then
{\small
$$
Sc_N(Y)= -\tfrac{1}{32}\,{\tfrac {32-279\,{Y}^{5}+66\,{Y}^{2}+1044\,{Y}^{4}-984\,{Y}^{3}-66\,{Y}^{6}+{Y}^{9}-12\,{Y}^{8}+54\,{Y}^{7}-432\,Y}{ \left( Y-3 \right) ^{2}Y^{2}}},
$$
}

observe that $Sc_N'(1)=0$ and $Sc_N''(1)=0$, but $Sc_N'''(1)=-\tfrac{9}{2}\neq 0$. Therefore, $y=1$ is an inflection point. 
\end{proof}

\section{Einstein Metrics on \texorpdfstring{$F(5)$}{F(5)}}
In this section, we will study the Eintein metrics on the full flag $F(5)=SU(5)/T^4$, in this particular case $\dim(\mathcal{M}_1^G)=9$. 

\subsection*{Standard metric}

Now we are going to study the standard metric,
a simple calculation shows that $2\rho=\tfrac{7}{10}$ and $Sc_N=7$. The
spectrum of $[L_\pg]$ is $\{0,\tfrac{1}{2},\tfrac{4}{5}\}$, with multiplicities 1,4,5 respectively. Therefore the standard metric is non-degenerate and unstable.

We observe that $0,\,\tfrac{1}{2}\,<2\rho\,<\,\tfrac{4}{5}$, so the coindex is 4.
We resume the results on the following proposition.
\begin{proposition}
 The standard metric on $F(5)$ is non-degenerate, unstable and its coindex is 4.   
\end{proposition}
\begin{remark}
    A more general proof of the last proposition can be found in \cite[Proposition 6.1]{jLauret}.
\end{remark}
\subsection*{Arvanitoyeorgos metric}
The Arvanitoyeorgos metric $g_A$ is given as
$$
(1,1,1,1,\tfrac{3}{2},\tfrac{3}{2},\tfrac{3}{2},\tfrac{3}{2},\tfrac{3}{2},\tfrac{3}{2}),
$$
up to scaling (see \eqref{classic-etrics}).
Consequently, we obtain $2\rho=\tfrac{11}{80}$ and $Sc_N=\tfrac{11}{4}\,{3}^{\frac{3}{5}}{2}^{\frac{2}{5}}\approx 7.0148$. 
It is important to note that the stability of this metric has not been addressed in previous studies.
{\tiny
$$
[L_\pg]=\left(\begin {array}{cccccccccc} \tfrac{2}{5}&{\tfrac {1}{60}}&{\tfrac {1}{60}}
&{\tfrac {1}{60}}&-{\tfrac {3}{20}}&-{\tfrac {3}{20}}&-{\tfrac {3}{20}}&0&0
&0\\\noalign{\medskip}{\tfrac {1}{60}}&\tfrac{2}{5}&{\tfrac {1}{60}}&{\tfrac {1}{
60}}&-{\tfrac {3}{20}}&0&0&-{\tfrac {3}{20}}&-{\tfrac {3}{20}}&0
\\\noalign{\medskip}{\tfrac {1}{60}}&{\tfrac {1}{60}}&\tfrac{2}{5}&{\tfrac {1}{60}
}&0&-{\tfrac {3}{20}}&0&-{\tfrac {3}{20}}&0&-{\tfrac {3}{20}}
\\\noalign{\medskip}{\tfrac {1}{60}}&{\tfrac {1}{60}}&{\tfrac {1}{60}}& \tfrac{2}{5}&0&0&-{\tfrac {3}{20}}&0&-{\tfrac {3}{20}}&-{\tfrac {3}{20}}
\\\noalign{\medskip}-{\tfrac {3}{20}}&-{\tfrac {3}{20}}&0&0&{\tfrac {17}{
30}}&-\tfrac{1}{15}&-\tfrac{1}{15}&-\tfrac{1}{15}&-\tfrac{1}{15}&0\\\noalign{\medskip}-{\tfrac {3}{20}}&0&-
{\tfrac {3}{20}}&0&-\tfrac{1}{15}&{\tfrac {17}{30}}&-\tfrac{1}{15}&-\tfrac{1}{15}&0&-\tfrac{1}{15}
\\\noalign{\medskip}-{\tfrac {3}{20}}&0&0&-{\tfrac {3}{20}}&-\tfrac{1}{15}&-\tfrac{1}{15}&
{\tfrac {17}{30}}&0&-\tfrac{1}{15}&-\tfrac{1}{15}\\\noalign{\medskip}0&-{\tfrac {3}{20}}&-
{\tfrac {3}{20}}&0&-\tfrac{1}{15}&-\tfrac{1}{15}&0&{\tfrac {17}{30}}&-\tfrac{1}{15}&-\tfrac{1}{15}
\\\noalign{\medskip}0&-{\tfrac {3}{20}}&0&-{\tfrac {3}{20}}&-\tfrac{1}{15}&0&-\tfrac{1}{15}&-\tfrac{1}{15}&{\tfrac {17}{30}}&-\tfrac{1}{15}\\\noalign{\medskip}0&0&-{\tfrac {3}{20}
}&-{\tfrac {3}{20}}&0&-\tfrac{1}{15}&-\tfrac{1}{15}&-\tfrac{1}{15}&-\tfrac{1}{15}&{\tfrac {17}{30}}
\end {array} \right).
$$
}

Whose eigenvalues are 
$$
\left\lbrace
0,\,\tfrac {19}{160}-{\tfrac {1}{480}}\sqrt {769},\,\tfrac {7}{40},\,\tfrac {19}{160}+{\tfrac {1}{480}}\,\sqrt {769},\,\tfrac{3}{16}
\right\rbrace,
$$
with multiplicities $1,\, 3,\, 2,\,3$ and $1$ respectively.  We can see that 
$$
\tfrac {19}{160}-{\tfrac {1}{480}}\sqrt {769}\,<2\rho\,<\tfrac{7}{40}\,<\tfrac {19}{160}+{\tfrac {1}{480}}\sqrt {769}\,<\tfrac{3}{16}.
$$
So $\lambda_p=\tfrac {19}{160}-{\tfrac {1}{480}}\sqrt {769}$, $\lambda_p^{\max}=\tfrac{3}{16}$ and the coindex is 3.
Finally, with the discussion above, we have proved the following.
\begin{proposition}
    The Arvanitoyeorgos metric $g_A$ is non-degenerate and unstable, and its coindex is 3.
\end{proposition}

\subsection*{Kähler metric \texorpdfstring{$g_K$}{gK}}
Let $g_K=(1,1,2,2,2,1,3,3,1,4)$ be the Kähler metric on $F(5)$, (see \eqref{classic-etrics}). A simple calculation shows that $2\rho=\tfrac{2}{5}$, and $Sc_N=4\,\sqrt {2}\sqrt [5]{3}\approx 7.0470. $

{\tiny
$$
[\Le_\pg(g_K)]= \left( \begin {array}{cccccccccc} {\tfrac {7}{30}}&\tfrac{1}{10}&-\tfrac{1}{5}&\tfrac{1}{15}&-\tfrac{1}{5}
&\tfrac{1}{10}&-\tfrac{1}{10}&0&0&0\\\noalign{\medskip}\tfrac{1}{10}&{\tfrac {7}{30}}&\tfrac{1}{15}&-\tfrac{1}{5}&-\tfrac{1}{5}&0&0&-\tfrac{1}{10}&\tfrac{1}{10}&0\\\noalign{\medskip}-\tfrac{1}{5}&\tfrac{1}{15}&{\tfrac {7}{12}}&\tfrac{1}{20}&0
&-\tfrac{1}{5}&0&-\tfrac{1}{5}&0&-\tfrac{1}{10}\\\noalign{\medskip}\tfrac{1}{15}&-\tfrac{1}{5}&\tfrac{1}{20}&{\tfrac {7}{12}}
&0&0&-\tfrac{1}{5}&0&-\tfrac{1}{5}&-\tfrac{1}{10}\\\noalign{\medskip}-\tfrac{1}{5}&-\tfrac{1}{5}&0&0&\tfrac{2}{3}&\tfrac{1}{15}&-\tfrac{1}{5}&
-\tfrac{1}{5}&\tfrac{1}{15}&0\\\noalign{\medskip}\tfrac{1}{10}&0&-\tfrac{1}{5}&0&\tfrac{1}{15}&{\tfrac {3}{20}}&\tfrac{1}{20}
&-\tfrac{1}{10}&0&-\tfrac{1}{15}\\\noalign{\medskip}-\tfrac{1}{10}&0&0&-\tfrac{1}{5}&-\tfrac{1}{5}&\tfrac{1}{20}&\tfrac{3}{4}&0&-\tfrac{1}{10}
&-\tfrac{1}{5}\\\noalign{\medskip}0&-\tfrac{1}{10}&-\tfrac{1}{5}&0&-\tfrac{1}{5}&-\tfrac{1}{10}&0&\tfrac{3}{4}&\tfrac{1}{20}&-\tfrac{1}{5}
\\\noalign{\medskip}0&\tfrac{1}{10}&0&-\tfrac{1}{5}&\tfrac{1}{15}&0&-\tfrac{1}{10}&\tfrac{1}{20}&{\tfrac {3}{20}}&-\tfrac{1}{15}\\\noalign{\medskip}0&0&-\tfrac{1}{10}&-\tfrac{1}{10}&0&-\tfrac{1}{15}&-\tfrac{1}{5}&-\tfrac{1}{5}&-\tfrac{1}{15}&{\tfrac 
{11}{15}}\end {array} \right).
$$
}

Whose spectrum is given by 
$$
\left\lbrace \lambda_\pg=0,\tfrac{47}{60}-\tfrac{\sqrt{85}}{60},\tfrac{47}{60}+\tfrac{\sqrt{85}}{60},\lambda_0,\lambda_1,\lambda_2,\lambda_\pg^{\text{max}}\right\rbrace,
$$
with multiplicities 4, 1, 1, 1, 1, respectively, where
$$
\begin{aligned}
5\lambda_p^{\text{max}}&={\tfrac {49}{12}}+{\tfrac {1}{72}}\,\sqrt {6}\sqrt {374+\sqrt [3]{
k_0}+57640\,{\tfrac {1}{\sqrt [3]{k_0}}}}\\
&+{\tfrac {1}{72}}\,\sqrt{4488-6\,\sqrt [3]{
k_0}-345840\,{\tfrac {1}{\sqrt [3]{k_0}}}+7056\,{\tfrac {\sqrt {6}}{\sqrt {374+\sqrt [
3]{k_0}+57640\,{\tfrac {1}{\sqrt [3]{
k_0}}}}}}},
\end{aligned}
$$
and
$$
k_0=7216460+420\,i\sqrt {790385991}
$$

A numerical analysis shows that
{\small
$$
0\,<\,2\rho=0.4<\lambda_0\,<\,\tfrac{47}{60}-\tfrac{\sqrt{85}}{60}\approx 0.6297\,<\,\lambda_1\,<\,\lambda_2\,<\,\tfrac{47}{60}+\tfrac{\sqrt{85}}{60}\approx 0.9370\,<\,\lambda_p^{\max}\approx 1.1496.
$$
}
So we have got the following proposition.
\begin{proposition}
The metric $g_K$ is non-degenerate and unstable, and its coindex is 3.
\end{proposition}
\subsubsection{Prescribed Ricci curvature problem}
The prescribed Ricci curvature problem consists of a given symmetric 2-tensor, asking about the existence and uniqueness (up to scaling) of a metric $g$ and a constant $c>0$ such that $\operatorname{Rc}(g)=cT$. In our case, this question is equivalent to asking about the injectivity (up to scaling) of $\operatorname{Rc}:\mca^G\to\mathcal{S}^2(M)^G$.  We say that a metric $g$ is Ricci locally invertible if $\operatorname{Rc}$ is invertible in an open neighborhood of $g$. According to \cite[Lemma 4.5.]{jLauret} $\left.d\operatorname{Rc}\right|_g=\tfrac{1}{2}\Le_\pg$, therefore a
necessary condition for the Ricci local invertibility of $g$ is $\dim \Ker \left.\Le_\pg\right|_{\mathcal{S}(M)^G}=1$.

The above analysis shows that $\lambda_p=0$, so the metric $g_K$ is not Ricci locally invertible, which motivates to search a curve of metrics $g(t)$, such that 

$$
Rc(g(t))=Rc(g_K)=\tfrac{2}{5}g_K=(\tfrac{2}{5},\tfrac{4}{5},\tfrac{4}{5},\tfrac{4}{5},\tfrac{2}{5},\tfrac{6}{5},\tfrac{6}{5},\tfrac{2}{5},\tfrac{8}{5}).
$$ 
A curve $g$ that satisfies this is defined by
$$
g(t)=(1-t,1+t,2-t,2+t,2,1,3,3,1,4),\qquad \qquad -1<t<1.
$$ 
\subsection*{Metrics \texorpdfstring{$g_1$}{g1} and \texorpdfstring{$g_2$}{g2}}
In this section, we provide two new Einstein metrics on $F(5)$ and study their stability properties.  We consider 
$$
g_i=(2C_i,2C_i,2C_i,2C_i,X_i,Y_i,X_i, X_i,Y_i,X_i), \qquad i=1,2, 
$$
where
$$
X_i=-\tfrac{10}{3}C_i^3+\tfrac{250}{9}C_i^2-\tfrac{2110}{27}C_i+\tfrac{2258}{27},
\qquad
Y_i=\tfrac{20}{3}C_i^3-\tfrac{4516}{27}+\tfrac{4760}{27}C_i-\tfrac{536}{9}C_i^2,
$$
and $C_1,C_2$ are the only two positive real roots of the polynomial
$$
P(x)=27x^4-270x^3+1008x^2-1700x+1129,
$$
which are explicitly given by 
{\small $$
\begin{aligned}
C_2,C_3&=\pm\tfrac{1}{2}\sqrt{-
 2 \sqrt[3]{\tfrac{\sqrt{7689}}{972} + \tfrac{13}{108}} - \tfrac{10}{27 \sqrt[3]{\tfrac{\sqrt{7689}}{972} + \tfrac{13}{108}}} + \tfrac{2}{9} + \tfrac{70}{27 \sqrt{\tfrac{1}{9} + \tfrac{10}{27 \sqrt[3]{\tfrac{\sqrt{7689}}{972} + \tfrac{13}{108}}}
 + 2 \sqrt[3]{\tfrac{\sqrt{7689}}{972} + \tfrac{13}{108}}}}}\\
 &+ \tfrac{1}{2}\sqrt{\tfrac{1}{9} + \tfrac{10}{27 \sqrt[3]{\tfrac{\sqrt{7689}}{972} + \tfrac{13}{108}}} 
+ 2 \sqrt[3]{\tfrac{\sqrt{7689}}{972} + \tfrac{13}{108}}} + \tfrac{5}{2}.
\end{aligned}
$$
}
Numerical approximations of the coordinates of the metrics are respectively given by
$$
(2C_2,X_2,Y_2)\approx(6.9152,7.7224,5.8870), 
$$
and 
$$
(2C_3,X_3,Y_3)\approx(5.8586,9.2800,5.7028).
$$
Aproximations to the metrics have been predicted in \cite{linoGrama}.
We provide the following approximate values for the normalized scalar curvatures:
$$
\scalar_N(g_1)\approx 7.0041, \qquad \scalar_N(g_2)\approx 6.9984.  
$$

\begin{proposition}\label{g2g3-stab}
The metrics $g_1$ and $g_2$ are both Einstein with Einstein constant $\rho=\tfrac{1}{10}$, they are non-degenerate and unstable with coindex $6$ and $5$, respectively.   
\end{proposition}

\begin{proof}
With the help of any computer algebra system, it is straightforward to check using \eqref{eq:Einstein Eq} that these metrics are Einstein with $\rho=\tfrac{1}{10}$, and using \eqref{prop: eq for F(n)} that the extremal eigenvalues of $[\Le_\pg]$, are $\lambda_\pg\approx 0.2695$ and $\lambda_\pg^{\text{max}}\approx 0.7422$.

A numerical analysis shows that in both cases, $\text{mult}(\lambda_\pg)=1$, $\text{mult}(\lambda_\pg^{\text{max}})=2$, $2\rho=\tfrac{1}{5}$, $\lambda_\pg<\tfrac{1}{5}<\lambda_\pg^{\text{max}}$ and the coindex of $g_1$ is $6$ and the coindex of $g_2$ is $5$.
\end{proof}

\subsection*{Metric \texorpdfstring{$g_3$}{g3}}
In this section, we present the following new Einstein metric on the manifold $F(5)$ and examine its stability. 

$$
g_3=(-\tfrac{3967}{5}+237K-\tfrac{87}{5}K^2,K,K,K,K,K,K,Z,Z,Z),
$$
where
$$
Z=\tfrac{41}{5}K^2-\tfrac{227}{2}K+\tfrac{3967}{10},
$$ and $K$ is the positive real root of $P(x)=14x^3-277x^2+1819x-3967$, explicitly, we have
$$
K=\tfrac{331}{1764 \sqrt[3]{\tfrac{5 \sqrt{1473}}{3528} + \tfrac{3623}{37044}}} + \sqrt[3]{\tfrac{5 \sqrt{1473}}{3528} + \tfrac{3623}{37044}} + \tfrac{277}{42}.
$$
Approximate values of the metric are:
{
$$
(-\tfrac{3967}{5}+237K-\tfrac{87}{5}K^2,K,Z)\approx(5.8091,7.4806,6.5187),
$$
}
and $Sc_N\approx 6.9988$.

\begin{proposition}
The metric $g_3$ is an Einstein metric with $\rho=\tfrac{1}{10}$, and it is non-degenerate, unstable and its coindex is 5. 
\end{proposition}
\begin{proof}
    The eigenvalues of $[\operatorname{L}_\pg]$ have a complicated expression, so we will provide only a numerical approximation of $\lambda_p$ and $\lambda_p^{\text{max}}$:
    
    $$
        \lambda_p(g_3) \approx 0.0557 \qquad \lambda_p^{\text{max}}(g_3)\approx 0.1389.
    $$
    
    To finish the proof, a numerical analysis shows that $\text{mult}(\lambda_p)=1$ and $\text{mult}(\lambda_p^{\text{max}})=2$, and $\lambda_p<\tfrac{1}{10}<\lambda_p^{\text{max}}$. 
    So the metric is non-degenerate and unstable. The analysis of the coindex is similar.
\end{proof}

\subsection*{Metrics \texorpdfstring{$g_4$}{g4} and \texorpdfstring{$g_5$}{g5}}
We consider in this section the following metrics:
$$
g_{i}=(M_i,M_i,\tfrac{1}{2}c_i,N_i,M_i,\tfrac{1}{2}c_i,N_i,\tfrac{1}{2}c_i,N_i,L_i) \quad i=4,5,
$$
where
{\small
$$
M_i=\tfrac{8971375}{1666862}c_i^5-\tfrac{98032375}{3333724}c_i^4+\tfrac{446084445}{6667448}c_i^3-\tfrac{68439950}{833431}c_i^2+\tfrac{357123741}{6667448}c_i-\tfrac{109792137}{8334310},
$$
}
{\small
$$
N_i=\tfrac{268587215}{13334896}-\tfrac{930186785}{13334896}c_i+\tfrac{1387293825}{13334896}c_i^2-\tfrac{1110422865}{13334896}c_i^3+\tfrac{239341375}{6667448}c_i^4-\tfrac{21499625}{3333724}c_i^5,
$$
}
and
{\small
$$
L_i=\tfrac{10670625}{3333724}c_i^5-\tfrac{129829875}{6667448}c_i^4+\tfrac{654761925}{13334896}c_i^3-\tfrac{876763875}{13334896}c_i^2+\tfrac{627815565}{13334896}c_i-\tfrac{183414881}{13334896}.
$$
}
The constant $c_i$ is a root of 
$$P(x)=17500x^6-113750x^5+316325x^4-490500x^3+447800x^2-225550x+48503.$$
In this case we have only got two real roots, $c_4$ and $c_5$.
The metrics can be approximated as follows:
{
$$
(M_4,\tfrac{1}{2}c_4,N_4,L_4)\approx(0.6470617,0.653643651,0.9237569,0.4766133),
$$
}
and
{
$$
(M_5,\tfrac{1}{2}c_5,N_5,L_5)\approx(0.6470616,0.923756890,0.653643,0.4766133).
$$
}
\begin{remark}
     Actually the metric $g_4$ and $g_5$ are isometric.   
\end{remark}

Although the previous observation allows us to focus solely on $g_4$.

$$
Sc_N = 7.0044
$$

\begin{proposition}
    The metric $g_4$ is an Einstein metric with Einstein constant $2\rho=\tfrac{1}{10}$. The metric is non-degenerate and unstable and its coindex is $4$.
\end{proposition}
\begin{proof}
    Using a mathematical symbolic software we see that $2\rho=\tfrac{1}{10}$. To see the $G$-stability, the spectrum of $[L_p(g_4)]$ is $\lbrace 0,\,\lambda_1,\,\lambda_2,\,\lambda_3,\,\mu_1,\,\mu_2,\,\mu_3\rbrace$, with $\text{mult}(\lambda_i)=1$  and $\text{mult}(\mu_i)=2$. 
    
    Finally a numerical analysis shows that the metric is non-degenerate, unstable and its coindex is $4$.
    
\end{proof}

\subsection{Non-homothetic metrics}

The aim of this section is to prove the following result.
\begin{theorem}
    The standard metric, the Arvanitoyeorgos metric, the K\"ahler metric and $g_1,\ldots,g_4$ are pairwise non-homothetic.
\end{theorem}
\begin{proof}    
    If we see the coindex of the metrics in Table \eqref{tablaMetricas}, then we only need to prove that 
    \begin{itemize}
        \item the standard metric is non-homothetic to $g_4$.
        \item Arvanitoyeogos metric is non-homothetic to $g_K$.
        \item $g_2$ is non-homothetic to $g_3$.
    \end{itemize}

    For the first case, we have that $[\Le_\pg](g_4)$ has an eigenvalue with multiplicity 2, and for the standard metric $[\Le_\pg](S)$ does not have an eigenvalue with such multiplicity, so they are non-homothetic.\\
    For the second case, we have that the eigenvalues 0 of $[\Le_\pg](g_K)$ verifies that $\text{mult}(0)>1$, and for the Arvanitoyeogos metric $\text{mult}(0)=1$, so they are non-homothetic.\\
    For the last case $g_2$ and $g_3$, a numeric analysis shows that 
    $$
    \tfrac{\lambda_\pg^{\max}(g_2)}{\lambda_\pg(g_2)}=6.016250172\neq 2.494623553=\tfrac{\lambda_\pg^{\max}(g_3)}{\lambda_\pg(g_3)},
    $$ 
    so those metrics are non-homothetic.     
\end{proof}

\begin{remark}
    Surprisingly, the approximations for $Sc_N$ of the seven metrics all fall within a range of less than 0.05. When two metrics have similar approximations for $Sc_N$,  it is generally difficult to determine if the metrics are not homothetic. In this case, we used an alternative and simple method to rule out the isometric case.
\end{remark}

\begin{table}
    \centering
\begin{tabular}{| c | c | c | c | c | c |}
\hline
Name  & Propieties & $2\rho$ & $Sc_N$ & coindex\\ \hline
$g_s$  & Non-deg. and unstable. & {\tiny$\tfrac{7}{10}$} & 7 & 4\\ \hline
$g_A$  & Non-deg. and unstable. & {\tiny$\tfrac{11}{80}$} & 7.0148 & 3\\ \hline
$g_K$  & Non-deg. and unstable mult(0)=4 & {\tiny$\tfrac{2}{5}$} & 7.0469 & 3\\ \hline
$g_1$  & Non-deg. and unstable & {\tiny$\tfrac{1}{10}$} & 7.0041 & 6\\ \hline
$g_2$  & Non-deg. and unstable & {\tiny$\tfrac{1}{10}$} & 6.9985 & 5\\ \hline
$g_3$  & Non-deg. and unstable & {\tiny$\tfrac{1}{10}$} & 6.9988 & 5\\ \hline
$g_4$  & Non-deg. and unstable & {\small$1$} & 7.0044 & 4\\ \hline
%$g_5$ &$(M_2,M_2,\tfrac{1}{2}c_2,N_2,M_2,\tfrac{1}{2}c_2,N_2,\tfrac{1}{2}c_2,N_2,L_2)$ & No deg. Unstable & 1 & 7.0043839 & 4\\ \hline
\end{tabular}
    \caption{Table with metrics on $F(5)$.}
    \label{tablaMetricas}
\end{table}

\bibliography{bibliography.bib} 
\bibliographystyle{amsplain}
\end{document}